\documentclass{amsart}
\usepackage{amsmath, amsthm, amssymb, mathtools}
\usepackage{algorithm,algorithmicx,algpseudocode}
\usepackage{graphics}
\usepackage{enumerate}
\usepackage{accents}
\usepackage{caption}
\usepackage{comment}
\usepackage{bm}
\usepackage{breqn}
\usepackage{cases}
\usepackage{eucal}
\usepackage{graphics}
\usepackage{fancyvrb}
\usepackage{seqsplit}
\usepackage{autobreak}
\usepackage{color}
\usepackage{ulem}
\usepackage{url}
\usepackage{booktabs} 

\usepackage[
        hypertexnames=false,
        hyperindex,
        pagebackref,
        breaklinks=true,
        bookmarks=false,
        colorlinks,
        linkcolor=blue,
        citecolor=red,
        urlcolor=red,
]{hyperref}

\theoremstyle{plain}
\newtheorem{Thm}{Theorem}[section]
\newtheorem{Prop}[Thm]{Proposition}
\newtheorem{Lem}[Thm]{Lemma}
\newtheorem{Cor}[Thm]{Corollary}
\newtheorem{Main}{Main Theorem}
\newtheorem{Main2}{Main Theorem}

\theoremstyle{definition}
\newtheorem{Def}[Thm]{Definition}
\newtheorem{Rem}[Thm]{Remark}
\newtheorem{Alg}[Thm]{Algorithm}
\newtheorem{Conj}[Thm]{Conjecture}
\newtheorem{Expect}[Thm]{Expectation}
\newtheorem{Proc}[Thm]{Procedure}

\makeatletter

\let\c@algorithm\c@Thm
\let\cl@algorithm\cl@Thm

\makeatother

\newcommand{\bA}{\mathbb{A}}
\newcommand{\bP}{\mathbb{P}}
\newcommand{\hS}{\mathcal{S}}

\newcommand{\Jac}{{\rm Jac}}

\DeclareRobustCommand{\erase}{\bgroup\markoverwith{{\rule[.5ex]{2pt}{0.4pt}}}\ULon}

\title{Superspecial genus-$4$ double covers of elliptic curves}
\author{Takumi Ogasawara}
\thanks{Yokohama National University, E-mail address: \texttt{ogasawara-takumi-np@ynu.jp}}
\author{Ryo Ohashi}
\thanks{The University of Tokyo, E-mail address: \texttt{ryo-ohashi@g.ecc.u-tokyo.ac.jp}}
\author{Kosuke Sakata}
\thanks{The University of Tokyo, E-mail address: \texttt{sakata-kosuke-rb@g.ecc.u-tokyo.ac.jp}}
\author{Shushi Harashita}
\thanks{Yokohama National University, E-mail address: \texttt{harasita@ynu.ac.jp}}

\numberwithin{equation}{subsection}

\begin{document}
\begin{abstract}
In this paper we study genus-$4$ curves
obtained as double covers of elliptic curves.
Firstly we shall give explicit defining equations of such curves with explicit criterion for whether it is nonsingular, and show the irreducibility of the long polynomial determining whether the genus-4 curve is nonsingular or not, in any characteristic $\ne 2,3$.
Secondly, as an application, we enumerate superspecial genus-$4$ double covers of elliptic curves in small characteristic.

\end{abstract}
\maketitle

\section{Introduction}
In this paper, let $K$ be a field of characteristic $p > 0$ and $\overline K$ its algebraic closure.
A curve $C$ over $K$ is called {\it superspecial} if its Jacobian variety is isomorphic (without a polarization) over $\overline K$ to a product of supersingular elliptic curves.  
It is well-known that $C$ is superspecial
if and only if the Frobenius operator on {$H^1(C,\mathcal O)$} 
(whose matrix is called Hasse-Witt matrix)
is zero. From this fact, one can say that superspecial curves are ``the most special curves". As in the case of abelian varieties,
superspecial {curves} play very important roles in the study of the moduli space of {curves}, it {is} very important to study superspecial curves. Furthermore, in many cases
superspecial curves have many rational points compared with their genus
or have big automorphism groups.
Based on these facts, constructing these curves has yielded fruitful applications in coding theory and cryptography.

Let us review the history of proving the existence (and the enumeration)
of the superspecial curves {of genus $g$}. Deuring \cite{Deuring} proved
the existence of a supersingular elliptic curve
(also see Igusa \cite{Igusa}).
The existence for genus $g=2$ and $p\ge 5$
is due to Serre \cite{Serre1983} and Ibukiyama-Katsura-Oort \cite[Proposition 3.1]{IKO}.
We find a proof of the existence for $g=3$ and $p\ge 3$ in Oort \cite[5.12]{Oort-Hyp} and Ibukiyama \cite[Theorem 1]{Ibukiyama}.
Note that Oort's proof is geometrical, whereas Ibukiyama's one is {arithmetic}.

The existence problem of a superspecial curve in the genus-$4$ case is still open in general, although many years have passed since the genus-$3$ case was settled.
However, there are many partial results known.
There are several kinds of genus-$4$ curves,
characterized by gonality (hyperelliptic and non-hyperelliptic)
and also by types of automorphism groups.
In the hyperelliptic case, computational results are found in \cite{KH18}, \cite{OKH23}, \cite{KNT}, and \cite{OK23}.
It is known that any non-hyperelliptic curve $C$ of genus $4$ is the complete intersection of a quadric surface and a cubic surface in $\mathbb P^3$.
The case was discussed in \cite{KH16} and \cite{KH17}.
A Howe curve (of genus $4$) is a curve of genus $4$ obtained as a fiber product of two elliptic curves over ${\mathbb P}^1$ (cf.\ \cite[Section 2]{KHS}). Note that any Howe curve is non-hyperelliptic (cf.\ \cite[Lemma 2.1]{KHH}).
A computational result on the enumeration of superspecial Howe curves
(of genus $4$) is found in \cite{KHH}.
Note that Howe curves were used in
Kudo, Harashita and Senda \cite{KHS} to prove the existence of supersingular curves of genus 4 in an arbitrary characteristic.

In this paper, we deal with a genus-$4$ curve $C$
obtained as a double cover of {a genus-$1$} curve $E$.
{By Riemann-Hurwitz's formula, the cover $C \rightarrow E$ is branched at 6 points.}
Note that these curves make a wider class than that of Howe curves.
{Over $\overline{K}$, such a curve} $C$ is characterized by the property that
automorphism group of $C$ contains a cyclic subgroup $G$ of order $2$
and the quotient of $C$ by $G$ is {a genus-$1$} curve.
{Moreover we shall see in Proposition \ref{prop:RealizationOfDCEC} that such a $C$ over $\overline{K}$ is described as in \eqref{eq:DCECinIntro} below.
In particular, such a $C$ is non-hyperelliptic (Corollary \ref{cor:DCECisNonhyper}).}

Let $K$ be a field of characteristic $p \ne 2,3$.
{We consider a curve over $K$ defined as}
\begin{equation}\label{eq:DCECinIntro}
    C=V(P,Q)
\end{equation}
in $\bP^3={\rm Proj}\,K[x,y,z,w]$ with
\begin{align*}
P(x,y,z,w) &:= y^2 z-f(x,z),\nonumber\\
Q(x,y,z,w) &:= w^2-q(x,y,z),
\end{align*}
where 
\begin{align*}
f(x,z) &:= x^3+Axz + Bz^3 \in K[x,z],\\
q(x,y,z) &:= a_0 x^2 + a_1 x y + a_2 x z + a_3 y^2 + a_4 y z + a_5 z^2 \in K[x,y,z],
\end{align*}
{
We call such a curve a {\it DCEC} over $K$,
if it is nonsingular.
This is an abbreviation of a double cover of an elliptic curve.}

The first main theorem gives a necessary and sufficient condition for the curve $C$ to be nonsingular:
\begin{Main}\label{MainTheorem1}
The curve $C$ {as in \eqref{eq:DCECinIntro}} is nonsingular if and only if
\begin{enumerate}
    \item[\rm (i)] $f_1:=4A^3+27B^2 \ne 0$, and 
    \item[\rm (ii)] $f_2:= \Delta/(f_3)^2 \ne 0$, 
\end{enumerate}
where $\Delta$ is the discriminant of 
the resultant $G(x)$ of $P|_{z=1}$ and $q|_{z=1}$ considered as polynomials in $y$
and $f_3$ is a short factor of $\Delta$ defined by
\[
    f_3:= Aa_1^2 a_3 a_4  - Ba_1^3 a_3   - a_0 a_1 a_4^2  - a_1^3 a_5  + a_1^2 a_2 a_4  + a_3 a_4^3.
\]
Moreover, $f_2$ is an irreducible polynomial over arbitrary field $K$ of characteristic $0$ or $p\ge 5$.
\end{Main}
\noindent We remark that $f_2$ is a very long polynomial (consisting of 3381 terms in characteristic $0$) in $A,B,a_0,\ldots,a_5$, see \cite{webpage} for the explicit form of $f_2$.
{See Lemma \ref{lem:meaning_f_3} for the geometric meaning of $f_3$.}

{The latter assertion in Main Theorem \ref{MainTheorem1}} is proved with the following geometric approach:
After we considering $A,B$ as constants satisfying $f_1\ne 0$, 
we will show that the locus $V(f_2)$ 
in the space of above $C$'s
where $C$ is singular is geometrically irreducible, {in Section 3, using a moduli theoretic argument}.

As an application, we discuss the existence/enumeration of superspecial curves of the type above.
The dimension of the moduli space of our curves $C$
is $6$, which is equal to the dimension of the moduli space of curves of genus $3$. This fact leads us to expect
that there may be a theoretical approach to show
the existence of superspecial curve of this type (cf.\ Section \ref{sec:conclusions}),
but in this paper we restrict ourselves to a computational approach.
We discuss the case that the genus-$4$ curve 
has a unique elliptic quotient of degree two.
We abbreviate such curves as UDCEC (a unique double cover of an elliptic curve).
Otherwise (i.e. if a genus-$4$ curve has two  elliptic quotient of degree two, then), it would be better to treat it in a different way, see \cite{KHH} for the case of Howe curves.

The aim of this paper is to give an algorithm
enumerating isomorphism classes of superspecial UDCEC's.
The algorithm consists of the following three subalgorithms:
\begin{itemize}
\item {an} algorithm to list superspecial DCEC's (cf.\ Subsection \ref{sebsec:Enumeration_SSP-DCEC});
\item {an} algorithm to select UDCEC's among them
(cf.\ Subsection \ref{subsec:UDCEC});
\item {an} isomorphism-{testing} algorithm for (U)DCEC's (cf.\ Subsections \ref{subsec:WholeAlgorihm} and  \ref{subsec:isomtest}).
\end{itemize}
Our implementation of these algorithms proves the following theorem for $p\le 23$.
The bound of $p$ can be increased by improving 
the algorithms, the implementation and the specification of the machine and so on.
\begin{Main}\label{MainTheorem2}
In the following, we abbreviate “superspecial" as “s.sp.".
\begin{enumerate}
    \item[\rm (1)] There is no s.sp.\ UDCEC in characteristics $p$ for $p<11$.
    \item[\rm (2)] There {is} 1 isomorphism {class} of s.sp. UDCEC's in characteristic 11.
    \item[\rm (3)] There are 2 isomorphism classes of s.sp.\,UDCEC's in characteristic 13.
    \item[\rm (4)] There are 13 isomorphism classes of s.sp.\,UDCEC's in characteristic 17.
    \item[\rm (5)] There are 16 isomorphism classes of s.sp.\,UDCEC's in characteristic 19.
    \item[\rm (6)] There are 53 isomorphism classes of s.sp.\,UDCEC's in characteristic 23.
\end{enumerate}
See Section \ref{sec:computational_results} for examples of the defining equations of these curves.
In Section \ref{sec:conclusions}, we propose an expectation on the asymptotic behavior of these numbers of isomorphism classes of s.sp.\,UDCEC's,
see Expectation \ref{expect:AsymptoticBehavior}.
\end{Main}

This paper is organized as follows.
In Section \ref{sec:Preliminaries}, with review of the classification of genus-4 curves, we study the genus-4 curves obtained as double covers of elliptic curves {and a way to compute their Hasse-Witt matrices. We also}
give an algorithm determining whether the curve
has a unique elliptic quotient.
{In Section \ref{sec:ProofMainTheorem1}, we prove Main Theorem \ref{MainTheorem1}.}
Section \ref{sec:ProofMainTheorem2} is devoted to the {proof of Main Theorem \ref{MainTheorem2}. In Section \ref{sec:conclusions}, we give a {conclusion} together with Conjecture \ref{radical} and Expectation \ref{expect:AsymptoticBehavior}.}

\subsection*{Acknowledgements}
{We thank the reviewers for careful reading {an earlier version of this article} and helpful comments.}
This research was supported by JSPS Grant-in-Aid for Scientific Research (C) 21K03159. This research was also conducted under a contract of “Research and development on new generation cryptography for secure wireless communication services" among “Research and Development for Expansion of Radio Wave Resources (JPJ000254)", which was supported by the Ministry of Internal Affairs and Communications, Japan.
\section{Preliminaries}\label{sec:Preliminaries}
In this section, we review a classification of non-hyperelliptic curves of genus 4 and we study a description of a DCEC.
We also give a criterion for
whether a DCEC $C$ is a UDCEC (a DCEC which has a unique elliptic quotient of degree $2$) and give an algorithm determining whether two UDCECs are isomorphic.

\subsection{Non-hyperelliptic curves of genus four}
Let $K$ be an algebraically closed field of characteristic 0 or $p \geq 3$, and let $C$ be a non-hyperelliptic curve of genus $4$ defined over $K$.
Then, the canonical divisor on $C$ is very ample, that is, this divisor defines an embedding
\[
    C \rightarrow \mathbb{P}^3 = {\rm Proj}\,K[x,y,z,w].
\]
It is known that $C$ is defined by an irreducible cubic form $P$ and an irreducible quadratic form $Q$ in $x,y,z$, and $w$ (cf. \cite[Chapter\,IV, Example 5.2.2]{Har}). In addition, we may assume that any coefficient of $P$ and $Q$ belongs to $K$ by \cite[Section 2.1]{KH16}. Hence, the curve $C$ is written as
    $V(P,Q)$ with $P,Q \in K[x,y,z,w]$ with $\deg{P} = 3$ and $\deg{Q} = 2$.

There is another construction of non-hyperelliptic curves of genus 4, studied by Howe \cite[Section 2]{Howe}, see \cite[Section 2]{KHS} for the following definition:
\begin{Def}\label{Howe}
A {\it Howe curve} is defined to be the desingularization of the fiber product $E_1 \hspace{-0.5mm}\times_{\mathbb{P}^1}\hspace{-0.5mm} E_2$ over the projective line $\mathbb{P}^1$, where $E_1$ and $E_2$ are elliptic curves which share {exactly} one ramified points.
\end{Def}
Howe curves make a subclass of non-hyperelliptic curves of genus $4$.
Here is a characterization of a Howe curve by the automorphism group:
\begin{Prop}\label{HoweV4}
A non-hyperelliptic curve $C$ of genus 4 is isomorphic to a Howe curve if and only if a subgroup $G$ of the automorphism group of $C$ is isomorphic to the Klein-4 group $V_4$ and the genus of $C/H$ is $1$ for an order-two subgroup $H$ of $G$.
\end{Prop}
\begin{proof}
The ``only if"-part is obvious. To prove the ``if"-part,
we apply \cite[Theorem B]{KR} to $G$
and the three order-two subgroups $H_1,H_2,H_3$
of $G$; then we have
\begin{equation}\label{eq:Kani-Rosen}
\Jac(C)^2\times \Jac(C/G)^4 \sim \Jac(C/H_1)^2\times \Jac(C/H_2)^2\times \Jac(C/H_3)^2.
\end{equation}
Since $C$ is non-hyperelliptic,
the genus of $C/H_i$ is either of $1$ and $2$,
and one of them is $1$ by the assumption.
By \eqref{eq:Kani-Rosen}, the genus $C/G$ is $0$
and the multiple set of the genera of $C/H_i$ for $i=1,2,3$ is $\{1,1,2\}$. Without loss of generality, we may assume that the genera of $E_1:=C/H_1$ and $E_2:=C/H_2$ are $1$. Then there exists a morphism $C$ to {the desingularization of}
the fiber product of $E_1$ and $E_2$ over $C/G$.
This is an isomorphism, since there is no non-trivial subgroup of any order-two group.
\end{proof}

\subsection{Double covers of elliptic curves}\label{subsec:Preliminaries-DCEC}
Let $K$ be an algebraically closed field of characteristic $0$ or $\ge 3$.
{In this subsection, all polynomials and all curves are over $K=\overline{K}$.}
Consider the projective curve 
\begin{equation}\label{DCEC-subsec2.2}
C=V(P,Q)
\end{equation}
in $\bP^3={\rm Proj}\,K[x,y,z,w]$ with
\begin{align}
P(x,y,z,w) &:= y^2 z-f(x,z),\label{eq:P}\\
Q(x,y,z,w) &:= w^2 - q(x,y,z),\label{PQ}
\end{align}
where $f(x,z)$ is a homogeneous polynomial of degree $3$ with $x^3$-coefficient$=1$ and moreover $q(x,y,z)$ is a homogeneous quadratic polynomial.
{Let $E$ be the genus-1 curve defined by $V(P)$ in $\bP^2={\rm Proj}\,K[x,y,z]$, then there exists the morphism
\begin{equation}\label{double-cover}
    \phi: C \rightarrow E\,;\,(x:y:z:w) \mapsto (x:y:z)
\end{equation}
of degree $2$.
This is the quotient map by the involution $(x,y,z,w)\mapsto (x,y,z,-w)$.
One can see that 
if $C$ is non-singular, then the quotient curve $E$ is also non-singular.
}

{Up to isomorphism of $C$, we may write}
\begin{align*}
    f(x,z) &= x^3 + Axz^2 + Bz^3,\\
    q(x,y,z) &= a_0 x^2 + a_1 x y + a_2 x z + a_3 y^2 + a_4 y z + a_5 z^2,
\end{align*}
{for $A,B,a_0,\ldots,a_5 \in K$.} 
The following two statements describe the nonsingularity of $C$ in terms of $P$ and $q$.
\begin{Prop}\label{sing_tan}
Assume that $V(P)$ is nonsingular.
The curve $C$ is singular at $(x:y:z:w)=(x_0:y_0:z_0:w_0)$ if and only if $w_0=0$ and the tangent space of $V(P)$ at $(x_0:y_0:z_0)$ is contained in the tangent space of $V(q)$ at $(x_0:y_0:z_0)$.
\end{Prop}
\begin{proof}
The Jacobian matrix of $C$ is given as
\[
    \begin{pmatrix}
        -3x^2-Az^2 & 2yz & y^2-2Axz-3Bz^2 & 0\\
        q_x & q_y & q_z & 2w
    \end{pmatrix},
\]
where $q_x,q_y$ and $q_z$ denote the partial derivative of $q$ with respect to $x,y$ and $z$, respectively.
Note that $C$ is nonsingular if and only if the matrix is of rank $2$ {at every point of $C$}.


If $C$ is singular at $(x_0:y_0:z_0:w_0)$, then $w_0=0$, especially
$(x_0:y_0:z_0)$ belongs to $V(q)$.
The above Jacobian matrix is of rank $\le 1$
at $(x_0:y_0:z_0:0)$
if and only if $(q_x,q_y,q_z)$ at $(x_0:y_0:z_0)$
is the first row multiplied by a scalar.
This is equivalent to saying that the tangent space of $V(P)$ at $(x_0:y_0:z_0)$ is contained in
the tangent space of $V(q)$ at $(x_0:y_0:z_0)$.
\end{proof}

{
Here is a criterion for whether $C$ is nonsingular:
}
\begin{Cor}\label{distinct_6pt}
The curve $C$ is nonsingular if and only if
$V(P)$ is nonsingular and
$V(P)$ and $V(q)$ intersect at {6 distinct} points
in ${\mathbb P}^2={\rm Proj}\,K[x,y,z]$.
\end{Cor}
\begin{proof}
The Bézout's theorem tells us that two curves $E=V(P)$ and $D=V(q)$ intersect at $6$ points with multiplicity.
{The “if"-part would be well-known but we show it for the reader's convenience.}
We suppose that $C$ is singular. Then, there exists a point $p = (x_0 : y_0 : z_0)$ such that $T_p(E) \subset T_p(D)$ from Proposition \ref{sing_tan}, where $T_p(E)$ and $T_p(D)$ are tangent spaces of $E$ and $D$ at $p$, respectively. Since the multiplicity at $p$ of $E$ and $D$ is equal to or greater than $2$, we have $\#(E\cap D) \leq 5$.
The “only-if"-part also follows from Proposition \ref{sing_tan}.
\end{proof}

{We call $C$ a {\it DCEC} if $C$ is nonsingular.
We can also determine the genus of any DCEC.}

\begin{Cor}\label{genus4}
Assume that $C$ is non-singular. Then, the curve $C$ is of genus $4$.
\end{Cor}
\begin{proof}
We apply Riemann-Hurwitz's formula (cf. \cite[Chapter\hspace{1mm}II, Theorem\hspace{1mm}5.9]{Silverman}) to the morphism (\ref{double-cover}).
The ramification index of $\phi$ at a point $P$ on $C$
is $2$ if the value of $w$ at $P$ is zero,
and $1$ otherwise.
{Corollary \ref{distinct_6pt}}
tells us that $C$ has {exactly} $6$ points whose $w$-coordinate $=0$.
Then, Riemann-Hurwitz's formula
implies that $C$ is a genus-$4$ curve.
\end{proof}

The next proposition gives a classification of genus-$4$ DCEC's
in terms of elliptic curve and six points on it.

\begin{Prop}\label{ClassificationOfDCEC}
Giving the following are equivalent to each other:
\begin{enumerate}
\item[\rm (i)] A nonsingular DCEC $C = V(P,Q)$ as in \eqref{DCEC-subsec2.2}.
\item[\rm (ii)] An elliptic curve $E$
and a set $\{P_1, P_2, \ldots, P_6\}$ of {six distinct} points {of} $E$ such that $\sum_{i=1}^6 P_i = O$.
\end{enumerate}
The correspondence is given by $E=V(P)$ and $\{P_1, P_2, \ldots, P_6\} = V(P,q)$,
where $q$ is the quadratic form in $x,y,z$ defined by $Q=w^2 - q$.
\end{Prop}
\begin{proof}
For a nonsingular DCEC $C = V(P,Q)$ where $P,Q$ are given as (\ref{eq:P}) and (\ref{PQ}), we have six points $V(P,q)=\{P_1, P_2, \ldots, P_6\}$. Put $E=V(P)$. Recall from \cite[Chapter\hspace{1mm}III, Remark\hspace{1mm}3.5.1]{Silverman} the following exact sequence:
\begin{equation}\label{eq:Abel-Jacobi}
  1 \longrightarrow \overline{K^*} \longrightarrow \overline{K}(E)^* \xrightarrow{\textrm{div}} \textrm{Div}^0(E) \xrightarrow{\sigma} E \longrightarrow 0,
\end{equation}
where $\sigma$ sends $(R)-(O)$ to $R$.
Then, we have
\[
    P_1+P_2+\cdots+P_6=\sigma\Biggl(\hspace{1mm}\sum_{j=1}^{6}{(P_{i})} - 6(O)\Biggr)
    = (\sigma \circ \textrm{div})\hspace{-0.3mm}\left(\dfrac{q}{z^2}\right) = O.
\]

{Conversely}, we assume that an elliptic curve $E=V(P)$ and six points $P_1,P_2,\ldots,P_6$ on $E$ satisfying $\sum_{i=1}^6 P_i=O$ are given.
By the exact sequence \eqref{eq:Abel-Jacobi} again, there exists an element $h \in \overline{K}(E)^*$
such that $\textrm{div}(h)=\sum_{i=1}^{6}{(P_{i})} - 6(O)$.
Since $q:=z^2h$ is everywhere regular, then this defines an element of $H^0(E,{\mathcal O}_E(2))
\simeq H^0({\mathbb P}^2,{\mathcal O}_{{\mathbb P}^2}(2))$.
Hence $q$ is a quadratic form in $x,y$ and $z$.
Putting $Q := w^2-q$, we have a nonsingular DCEC $V(P,Q)$ of genus $4$.
\end{proof}


{Since $\deg{\phi} = 2$ where $\phi$ is defined as (\ref{double-cover}), a DCEC is
a double cover of a genus-$1$ curve. 
Conversely, over $K$ an algebraically closed field we have:
\begin{Prop}\label{prop:RealizationOfDCEC}
Let $C$ be a genus-4 double cover of a genus-$1$ curve. Then $C$
is a DCEC, i.e., is realized as in \eqref{DCEC-subsec2.2}.
\end{Prop}
\begin{proof}
Let $\phi : C \rightarrow E$ be a morphism to a genus-$1$ curve $E$ of degree $2$.
Since $\deg\phi =2$,
there exists $g \in K(E)$ such that $K(C) = K(E)(\sqrt{g})$.
The Riemann-Hurwitz formula tells us that the branch points (of ramification index $2$) of $C\to E$ are distinct $6$ points on $E$, say $P_1,\ldots,P_6$.
Hence ${\rm div}(g) = (P_1) + \cdots + (P_6) + 2 D$ for some divisor $D$ of $E$.
Choose a point $O'$ of $E$ and consider $E$ as an elliptic curve with origin $O'$.
Choose $O\in E$ so that $[6]O=P_1+\cdots + P_6$. Then
$(P_1) + \cdots + (P_6) - 6(O)$ is principal, and reconsider
$E$ as an elliptic curve with origin $O$.
We can write $g=u^2 h$ for some $u\in K(E)$, where ${\rm div}(h) = (P_1) + \cdots + (P_6) - 6(O)$
and ${\rm div}(u) = D + 3(O)$.
Then $C$ is isomorphic to the DCEC constructed from $P_1, \ldots, P_6$ on $E$ in
the proof of Proposition \ref{ClassificationOfDCEC}.
\end{proof}
}

Since any complete intersection of a quadratic and a cubic surfaces in ${\mathbb P}^3$ is canonical (\cite[Chap.\ IV, Example 5.2.2]{Har}), we have
\begin{Cor}\label{cor:DCECisNonhyper}
    Any {genus-4 double cover of a genus-$1$ curve} is non-hyperelliptic.
\end{Cor}

\subsection{Relation between Howe curve and DCEC}
In this subsection, we study a sufficient condition for a DCEC to be a Howe curve.
By definition one can easily see that a Howe curve is a {\rm DCEC}, but
this also can be checked by explicit defining equations. 
A Howe curve $C$ is {the desingularization of the curve defined by} the equations
\begin{equation}\label{eq:HoweCurve}
    \begin{split}
    y^2z - g_1(x,z) &= 0,\\
    w^2z - g_2(x,z) &= 0,
    \end{split}
\end{equation}
where $g_1$ and $g_2$ are coprime cubic polynomials with $x^3$-coefficient $1$. Then, we have that $h_1(x,z) := (g_1(x,z) - g_2(x,z))/z$ is a quadratic polynomial.
We put
\begin{align}
    P_1 &:= y^2z - g_1(x,z),\\
    Q_1 &:= w^2 + h_1(x,z) - y^2. \label{eq:HoweQ}
\end{align}
Then $V(P_1,Q_1)$ is isomorphic to $C$, whence
$C$ is a DCEC.
\begin{Prop}
If $C = V(P,Q)$ as defined in \eqref{DCEC-subsec2.2} is a nonsingular DCEC with $a_1 = a_4 = 0$ and $a_3\ne 0$, then $C = V(P,Q)$  is a Howe curve.
\end{Prop}
\begin{proof}
Let $C=V(P,Q)$ be a nonsingular DCEC as in \eqref{DCEC-subsec2.2} with $a_1=a_4=0$ and $a_3\ne 0$.
By replacing $Q$ as in \eqref{PQ} by $Q/a_3$ and $w$ by $w/\hspace{-0.3mm}\sqrt{a_3}$, then we may assume that $a_3=1$.
Then $Q$ is of the form of the right hand side of \eqref{eq:HoweQ}. Set 
\[
    h(x,z):=Q-w^2+y^2, \quad q:=w^2-Q.
\]
The polynomials $f$ and $h$ are coprime. Indeed, otherwise $P=0$ and $q=0$
have a double root at $y=0$, which contradicts that $C = V(P,Q)$ is nonsingular.
We put
\[
    g_1(x,z) := f, \quad g_2(x,z) := f-zh,
\]
then $g_1$ and $g_2$ are coprime, as $z$ is not a divisor of $f$.
Hence we have a Howe curve defined by \eqref{eq:HoweCurve}, which is isomorphic to $C$.
\end{proof}

\subsection{{Hasse-Witt} matrices and superspeciality}
Let $K$ be a perfect field of characteristic $p \geq 3$ and let $\overline{K}$ be the algebraic closure of $K$. 
Let $C = V(P,Q)$ be the nonsingular genus-4 curve in ${\mathbb P}^3={\rm Proj}\,K[x,y,z,w]$ defined by a cubic form $P$ and a quadratic form $Q$ in $x,y,z$ and $w$.
Put
{
\begin{equation}\label{eij}
(e_{ij}) = \begin{pmatrix}
1 & 1 & 1 & 2\\
1 & 1 & 2 & 1\\
1 & 2 & 0 & 1\\
2 & 1 & 1 & 1
\end{pmatrix}.
\end{equation}
}
Recall the following lemma on the {Hasse-Witt} matrix of $C$:
\begin{Lem}[{\cite[Proposition 3.1.4]{KH16}}]
Let $M$ be the {Hasse-Witt} matrix of $C$ with respect to a basis of $H^1(C,\mathcal O_C)$. Then, the $(i,j)$-entry of $M$ is the coefficient of
$x^{pe_{i1}-e_{j1}}y^{pe_{i2}-e_{j2}}z^{pe_{i3}-e_{j3}}w^{pe_{i4}-e_{j4}}$ in
$(PQ)^{p-1}$ for $1\le i,j \le 4$.
\end{Lem}

Since a nonsingular curve is superspecial if and only if its {Hasse-Witt} matrix is the zero matrix (cf. \cite[Thoerem 4.1]{Nygaard}), we have
\begin{Cor}\label{c-mMtx}
The curve $C$ is superspecial if and only if all the coefficients of $x^{pe_{i1}-e_{j1}}y^{pe_{i2}-e_{j2}}z^{pe_{i3}-e_{j3}}w^{pe_{i4}-e_{j4}}$ in $(PQ)^{p-1}$ is zero for $1\le i,j\le 4$ .
\end{Cor}
We have some strategies to reduce time complexity of computing {Hasse-Witt matrix} of $C$. The first idea is given by following lemma.

\begin{Lem}\label{s-sing}
Let $C=V(P,Q)$ be a curve given by (\ref{DCEC-subsec2.2}). If $C$ is nonsingular and superspecial, then the elliptic curve $E:=V(P)$ is supersingular.
\end{Lem}
\begin{proof}
This follows from \cite[Proposition 6]{Lachaud} and the fact
that any superspecial curve has an ${\mathbb F}_{p^2}$-form
which is maximal or minimal (cf. \cite[Theorem 1.1]{Ekedahl}).
\end{proof}
Actually, the $(1,1)$-entry of {Hasse-Witt matrix} of $C$ {, i.e., the $x^{p-1}y^{p-1}z^{p-1}w^{2p-2}$-coefficient of $(PQ)^{p-1}$ equals} the coefficient of $x^{p-1}$ in the polynomial $f(x,1)^{(p-1)/2}$ {up to sign}, where $f(x,z)$ is given by (\ref{eq:P}). 
Hence, computing supersingular elliptic curves in advance reduces the time complexity of computing the {Hasse-Witt} matrix of $C$. 

\begin{algorithm}[htbp]
\begin{algorithmic}[1]
\renewcommand{\algorithmicrequire}{\textbf{Input:}}
\renewcommand{\algorithmicensure}{\textbf{Output:}}
\algnewcommand{\IIf}[1]{\State\algorithmicif\ #1\ \algorithmicthen}
\algnewcommand{\EndIIf}{\unskip\ \algorithmicend\ \algorithmicif}
\caption{\ Computing {Hasse-Witt} matrix}\label{Cariter-Manin}
\label{carier-manin:algrm}
\Require A DCEC $C = V(P,Q)$, where $Q = w^2 - q(x,y,z)$
\Ensure The {Hasse-Witt} matrix of $C$
\For{$1 \leq i,j \leq 4$}
\State Compute $T_{ij} := pe_{i4} - e_{j4}$, where $e_{ij}$ is defined as (\ref{eij})
\If{$T_{ij}$ is even}
\State Compute $H := P^{p-1}\cdot S_{ij}$, where $S_{ij}$ is the $w^{T_{ij}}$-coefficient of $Q^{p-1}$
\[
    S_{ij} := \binom{p-1}{T_{ij}/2}(-q(x,y,z))^{p-1-T_{ij}/2}
\]
\State Set $c_{ij}$ be the $x^{pe_{i1}-e_{j1}}y^{pe_{i2}-e_{j2}}z^{pe_{i3}-e_{j3}}$-coefficient of $H$
\Else
\State Set $c_{ij} := 0$
\EndIf
\EndFor\\
\Return $(c_{ij})_{i,j}$
\end{algorithmic}
\end{algorithm}

Here are more ways to reducing
the complexity in computing the {Hasse-Witt} matrix of $C$. 
Firstly, we compute $P^{p-1}$ and $Q^{p-1}$ separately, instead of $(PQ)^{p-1}$. 
Secondly, we do not compute every term of $Q^{p-1}$. Since $P$ does not contain the indeterminate $w$, it suffices to compute only the term containing $w^{pe_{i4} - e_{j4}}$ of $Q^{p-1}$.
Summarizing the above discussions, we obtain Algorithm \ref{carier-manin:algrm} in order to compute the {Hasse-Witt} matrix of $C$.

\subsection{Criterion {for} whether DCEC is UDCEC}\label{subsec:UDCEC}
If a DCEC $C$ has a unique elliptic quotient of degree $2$, we call $C$ a {\it UDCEC}. Otherwise (i.e.,\ if a DCEC $C$ has more than one elliptic quotients of degree $2$), then $C$ is a Howe curve or another special curve. It would be better to consider such cases in different contexts, for example, the enumeration of superspecial Howe curves has been discussed in \cite{KHH}. Hence in this paper, we restrict ourselves to the case that $C$ has a unique elliptic quotient of degree two. Here is a criterion for whether $C$ has a unique elliptic quotient of degree two:
\begin{Lem}\label{uniqueornot}
The following are equivalent:
\begin{enumerate}
    \item[\rm (1)] A DCEC $C = V(P,Q)$ has a unique elliptic quotient $E = V(P)$.
    \item[\rm (2)] For any linear coordinate change $(x,y,z,w) \leftrightarrow (x',y',z',w')$, there does not exist $(\lambda,a,b,c,d)\in K^\times\times K^4$ with $(a,b,c,d)\ne (0,0,0,0)$ such that
    \[
        P' = \lambda P + (ax + by + cz + dw)Q
    \]
    is a polynomial only in $(x',y',z')$.
\end{enumerate}
\end{Lem}
\begin{proof}
If there is another elliptic quotient $C\to V(P')$ where $V(P')$ is an elliptic curve in ${\rm Proj}k[x',y',z']$, then we have
$V(P,Q)\simeq V(P',Q')$ for another form $V(P',Q')$ of the DCEC. Since DCEC is non-hyperelliptic (Corollary \ref{cor:DCECisNonhyper}) {and therefore $C$ is a canonical curve in ${\mathbb P}^3$},
the isomorphism is realized as a linear coordinate change.
After the linear coordinate change,
we may assume $V(P,Q) = V(P',Q') \subset {\mathbb P}^3$.
Since quadratic forms in the defining ideal is unique (up to multiplication by a nonzero constant), we may assume $Q=Q'$. Since any cubic form in the defining ideal 
has to be of the form $\lambda P + (ax+by+cz+dw)Q$, we have the lemma.
\end{proof}

Let $E' = V(P')$ be another elliptic quotient of $C = V(P,Q)$.
Consider the ring
\[
    R := K[\alpha_{ij},a,b,c,d\,;\,1\le i\le 3,\ 1\le j \le 4]
\]
and the equation
\begin{eqnarray}\label{Lem2-16}
    &&P'(\alpha_{11} x + \alpha_{12} y + \alpha_{13} z + \alpha_{14} w,
    \ldots, \alpha_{31} x + \alpha_{32} y + \alpha_{33} z + \alpha_{34} w)\\ \notag
    &&= P + (a x + b y + c z + d w)Q.
\end{eqnarray}
Let $I \subset R$ be the ideal generated by all the polynomials obtained by comparing the coefficients with respect to $x,y,z,w$ of both sides of (\ref{Lem2-16}). For each $s \in \{a,b,c,d\}$, define $I_s$ as the ideal of $R[t]$ generated by $I$ and $st-1$, where $t$ is a new parameter. If $V(I_s)$ is empty (i.e., the Gr\"obner basis of $I_s$ is equal to $\{1\}$) for all $s \in \{a,b,c,d\}$, then $C = V(P,Q)$ does not have $V(P')$ as another elliptic quotient of degree $2$.

\begin{algorithm}[htbp]
\begin{algorithmic}[1]
\renewcommand{\algorithmicrequire}{\textbf{Input:}}
\renewcommand{\algorithmicensure}{\textbf{Output:}}
\algnewcommand{\IIf}[1]{\State\algorithmicif\ #1\ \algorithmicthen}
\algnewcommand{\EndIIf}{\unskip\ \algorithmicend\ \algorithmicif}
\caption{Criterion for whether $C$ has a unique quotient} \label{checkUDCEC}
\Require A superspecial DCEC $C = V(P,Q)$ as in (\ref{DCEC-subsec2.2}).
\Ensure Whether $C$ is UDCEC or not
\State Let $\mathcal{E}$ be the list of the isomorphism classes of supersingular elliptic curves
\For{$E' \in \mathcal{E}\!\smallsetminus\{V(P)\}$}
\State Write $V(P') = E'$ with a polynomial $P' \in K[x,y,z]$
\State Let $L \subset K[a,b,c,d]$ be the list of $x^iy^jz^kw^l$-coefficients of
\[
    F := P + (ax+by+cz+dw)Q -P' \in K[a,b,c,d][x,y,z,w]
\]
\For{$s \in \{a,b,c,d\}$}
\State Let $I_s \subset K[a,b,c,d,t]$ be the ideal generated by the polynomials in $L$ and $st-1$~(where $t$ is a ``dummy" indeterminate).
\IIf{$I_s \ne \{1\}$} \ \,\Return \,False\ \ \EndIIf
\EndFor
\EndFor\\
\Return True
\end{algorithmic}
\end{algorithm}

Hence, it follows from {Lemma \ref{uniqueornot}} that we obtain Algorithm \ref{checkUDCEC} for checking whether a DCEC has a unique elliptic quotient.

\subsection{Isomorphism {testing} of UDCEC's}\label{subsec:isomtest}
Let $C_1 = V(P_1,Q_1)$ and $C_2 = V(P_2,Q_2)$ be UDCEC's as in (\ref{DCEC-subsec2.2}), and define $q_1 := w^2 - Q_1,\,q_2 := w^2 - Q_2$ as in (\ref{PQ}). In this subsection, we determine when $C_1$ and $C_2$ are isomorphic to each other.
\begin{Prop}\label{CriterionForIsom}
The UDCEC's $\,C_1$ and $C_2$ are isomorphic to each other if and only if
\begin{enumerate}
    \item[\rm (1)] the $j$-invariants of $E_1=V(P_1)$ and $E_2=V(P_2)$ are the same,
    \item[\rm (2)] there exist an isomorphism $\psi: E_1 \to E_2$ as genus-$1$ curves and $c \in \overline{K}\!\smallsetminus\!\{0\}$ such that $\psi^* q_2 = c q_1$. In case, with identifications of $E_1$
    and $E_2$ as elliptic curves, say $E:=E_1=E_2$, the isomorphism $\psi: E \to E$
    is the composition of a translation by an element of $E[6]$
    and an automorphism of $E$.
\end{enumerate}
\end{Prop}

\begin{proof}
We consider {the} six points $\{p_{i1},\ldots, p_{i6}\} \subset E_i$ {whose} sum is zero, {which determine $Q_i$ as in Proposition \ref{ClassificationOfDCEC}}.
If there exists an isomorphism $C_1 \rightarrow C_2$, then it induces an isomorphism $\psi: E_1\to E_2$ (since $C_1$ and $C_2$ are UDCEC's), whence $E_1$ and $E_2$ have the same $j$-invariant. Moreover, we have that 
\[
    \psi\left( \{p_{11},\cdots,p_{16}\}\right) = \{p_{21},\cdots,p_{26}\}.
\]
Conversely, the existence of such a $\psi$ implies that $C_1 \cong C_2$ (Proposition \ref{ClassificationOfDCEC}).

We {identify} $E_1$ and $E_2$ via an isomorphism as elliptic curves and write it as $E$.
Then there exist $\rho \in{\rm Aut}(E)$ and a point $R\in E$ such that $\psi = \tau_R \circ \rho$, where $\tau_R$ is the translation-by-$R$ map.
By $\rho = \tau_R^{-1} \circ \psi$ we have
\[
    \rho \Biggl(\hspace{1mm}\sum _{j=1}^6 {p_{1j}}\Biggr) = \sum _{j=1}^6{\rho ( {p_{1j}})} \\
= \sum _{j=1}^6{p_{2j}}  - 6R. 
\]
Since $\sum _{j=1}^6{p_{ij}}=O$, we have $R \in E[6]$.
\end{proof}

{Proposition \ref{CriterionForIsom} gives an algorithm for checking if two UDCEC's are isomorphic to each other or not}. Here, we write down an algorithm (Algorithm \ref{isomtest}) for this explicitly, for reader's {convenience}.
\begin{algorithm}[htbp]
\begin{algorithmic}[1]
\caption{Isomorphism {testing}} \label{isomtest}
\renewcommand{\algorithmicrequire}{\textbf{Input:}}
\renewcommand{\algorithmicensure}{\textbf{Output:}}
\algnewcommand{\IIf}[1]{\State\algorithmicif\ #1\ \algorithmicthen}
\algnewcommand{\EndIIf}{\unskip\ \algorithmicend\ \algorithmicif}
\Require Two UDCEC's $C_1 = V(P_1,Q_1)$ and $C_2 = V(P_2,Q_2)$ as in (\ref{DCEC-subsec2.2})
\Ensure Whether $C_1$ and $C_2$ are isomorphic or not
\State Let two elliptic curves $E_1 = V(P_1)$ and $E_2 = V(P_2)$
\If{$j(E_1) \neq j(E_2)$} \ \,\Return \,False
\Else
    \For{$P \in E_1[6]$ and $\phi \in {\rm Aut}(E_1)$}
        \State Let $\tau_P: E_1 \rightarrow E_1$ be the translation-by-$P$ map and $\psi := \tau_P \circ \phi$.
        \IIf{$\psi^*(q_2)/q_1 \in K^\times$ with $q_i = w^2 - Q_i$} \ \,\Return \,True\ \,\EndIIf
    \EndFor\\
    \Return \,False
\EndIf
\end{algorithmic}
\end{algorithm}

{The next corollary is important, since it tells us a field
over which all the elements of ${\rm Aut}(E_1)$ in Algorithm \ref{isomtest} are defined.}

\begin{Cor}
Let $C = V(P,Q)$ be a superspecial UDCEC {over $K$}. Then, the following are true:
\begin{enumerate}
    \item[\rm (1)] Let $E = V(P)$ be an elliptic curve over $K$. Then, any automorphism of $C$ 
    is defined over $K(\mu\!\!\!\!\mu_{12},E[6])$, where $\mu\!\!\!\!\mu_{12}$ is the group of $12$-th roots of unity.
    \item[\rm (2)] If $E$ is ${\mathbb F}_{p^2}$-minimal or ${\mathbb F}_{p^2}$-maximal with $K'={\mathbb F}_{p^2}$, then we have that
    \[
        K(\mu\!\!\!\!\mu_{12},E[6])=K\cdot{\mathbb F}_{p^4}.
    \]
\end{enumerate}
\end{Cor}
\begin{proof}
Write $E: y^2 = x^3 + Ax + B$ with $A,B \in K$. {It} is known that any automorphism of $E$ as elliptic curve is defined over $K(\mu\!\!\!\!\mu_{12})$. Note that $E$ is supersingular by Lemma \ref{s-sing}. Then, it follows from Proposition \ref{CriterionForIsom}\,(2) that any automorphism of $C = V(P,Q)$ is defined over $K(\mu\!\!\!\!\mu_{12},E[6])$. This completes the proof of (1).

If $E$ is ${\mathbb F}_{p^2}$-minimal or ${\mathbb F}_{p^2}$-maximal, then we have that $\phi^2= [\pm p]$ on $E$ {where $\phi$ is the Frobenius map on $E$}.
For any $P=(x_0,y_0) \in E[6]$, we have
\[
    \phi^2 P = \pm p P = P \text{ or }-P.
\]
Hence, we have $({x_0}^{p^2},{y_0}^{p^2}) = (x_0,\pm y_0)$.
This means that $P=(x_0,y_0)$ is an ${\mathbb F}_{p^4}$-rational point.
The second assertion of the corollary follows from 
$\mu\!\!\!\!\mu_{12} \subset {\mathbb F}_{p^4}$.
\end{proof}

\section{Proof of Main Theorem \ref{MainTheorem1}}\label{sec:ProofMainTheorem1}
In this section, we prove Main Theorem \ref{MainTheorem1} stated in Introduction. Let $K$ be an algebraically closed field of characteristic $p\ge 5$, and let $C$ be the projective curve 
\[
C=V(P,Q)
\]
in $\bP^3={\rm Proj}\,K[x,y,z,w]$ with
\begin{align*}
    P(x,y,z,w)&:= y^2 z-x^3-Axz^2-Bz^3,\\
    Q(x,y,z,w)&:= w^2-a_0 x^2-a_1 xy-a_2 xz-a_3 y^2  -a_4 yz-a_5 z^2
\end{align*}
for $A,B, a_0,\ldots,a_5 \in K$. Let $q(x,y,z) := a_0 x^2-a_1 xy-a_2 xz-a_3 y^2  -a_4 yz-a_5 z^2$.

\begin{Rem}\label{rem:p_NotTooSmall}
We assume $p\ge 5$ so that any elliptic curve is defined as $V(P)$ the above.
It would be possible to consider the case of $p=2,3$ by using {a} more general form for $P$,
but {with higher} computational costs.
\end{Rem}

In Corollary \ref{distinct_6pt}, we have seen that $C$ is nonsingular if and only if $V(P)$ is nonsingular, and $V(P)$ and $V(q)$ intersect at {6 distinct} points in ${\mathbb P}^2={\rm Proj}\,K[x,y,z]$.
The latter condition is equivalent to that $V(P)$ and $V(q)$ intersect transversely,
{
in other words the following two statements hold:
\begin{enumerate}
    \item $J(P,q) \ne 0$ at any point of $V(P,q)\setminus V(z)$. Here $J(P,q)$ is the determinant of the Jacobian matrix
    \[
        \begin{pmatrix}
        \dfrac{\partial P|_{z=1}}{\partial x}  & \dfrac{\partial P|_{z=1}}{\partial y}\\[2.5mm]
        \dfrac{\partial q|_{z=1}}{\partial x}  & \dfrac{\partial q|_{z=1}}{\partial y}\\  
        \end{pmatrix}.
    \]
    \item If $[0:1:0]\in V(P,q)$ or equivalently $a_3=0$, then the determinant of the Jacobian matrix at  $[0:1:0]$ is not zero.
\end{enumerate}
}

Before we look at an explicit condition for $C$ being singular,
we study the short factor $f_3$ of $\Delta$ in Main Theorem \ref{MainTheorem1}.
Before the next Lemma, we assume $a_3 \ne 0$.
We have $\Delta\ne 0$ if and only if
the $x$-coordinates of the points of $V(P|_{z=1},q|_{z=1})$ consists of distinct $6$ values,
in particular $C$ is nonsingular. But
$\Delta$ can be zero even if $C$ is nonsingular;
in case, for the $x$-coordinate $x_0$ of a point of $V(P|_{z=1},q|_{z=1})$ 
the simultaneous equation $P|_{z=1}(x_0,y)=0,q|_{z=1}(x_0,y)=0$ in $y$ has multiple roots (taking account of {multiplicity}),
equivalently $a_1x_0+a_4=0$. 
The next lemma says that this phenomenon occurs at the generic point of $V(f_3)$.
\begin{Lem}\label{lem:meaning_f_3}
We have
\[
    \left(P|_{z=1},q|_{z=1},a_1x+a_4\right) \cap K[A,B,a_0,\ldots,a_5] = (f_3).
\]
\end{Lem}
\begin{proof}
By a PC computation (\cite{webpage}),
the Gr\"obner basis of the ideal $(P|_{z=1},q|_{z=1},a_1 x + a_4, 6t-1)$ of
${\mathbb Z}[x,y,A,B,a_0,\ldots,a_5,t]$ with respect to the lexicographic order with
$x\succ y \succ \cdots \succ a_0 \succ \cdots \succ a_5  \succ t$
consist of $f_3$, $6t-1$ and polynomials with leading coefficient $1$ whose leading monomial contains $x$ or $y$.
This means that the Gr\"obner basis (without $6t-1$) is
the Gr\"obner basis $G$ of the ideal $(P|_{z=1},q|_{z=1},a_1 x + a_4)$ of
$K[x,y,A,B,a_0,\ldots,a_5]$ for any field $K$ of characteristic $\ne 2,3$.
The lemma follows from that $f_3$ is the unique member of $G$
containing neither $x$ nor $y$. 
\end{proof}

The next proposition gives a condition for $C$ being singular.

\begin{Prop}\label{prop:SingularLocus} Let $\infty$ denote the point {$[0:1:0:0]$ on ${\mathbb P}^3$}, and set 
\[
I := \left(P|_{z=1},q|_{z=1}, J(P,q)\right) \cap K[A,B,a_0,\ldots,a_5].
\]
Then we have 
\begin{enumerate}
    \item[\rm (1)] $V(I) =  V(f_2)$.
\item[\rm (2)]
    The curve $C$ is singular at $\infty$ if and only if $a_1=a_3=0$. Moreover $(f_2)\subset (a_1, a_3)$.
\end{enumerate}
Combining (1) and (2), we see that the locus in ${\rm Spec}(K[A,B,a_0,\ldots,a_5])$ where $C$ is singular
is a dense subset of $V(f_2)$.
\end{Prop}
\begin{proof}
(2) It is clear that $\infty = {[0:1:0:0]} \in C$ if and only if $a_3=0$.
Assume that $a_3=0$. By the Jacobian criterion, we see that $C$ is singular at $\infty$ if and only if $a_1 = 0$. It is straightforward to see $(f_2) \subset (a_1,a_2)$.

(1) 
By the computation on MAGMA~\cite{MAGMA} over $\mathbb Q$, we write $f_2$ in the form $r_1 P + r_2 q + r_3 J(P,q)$ for $r_i \in {\mathbb Q}[{x,y,} A,B,a_0,\ldots,a_5]$.
Specifically, we use the built-in functions {\tt IdealWithFixedBasis} and {\tt Coordinates} of MAGMA
for $[P|_{z=1},q|_{z=1},J(P,q)]$ and $f_2$,
and confirm that all the denominators of the coefficients of $r_i$
are products of $2$ and $3$, see \cite{webpage} for the code we used.
Hence
\begin{eqnarray}
    (f_2) \subset I,
\end{eqnarray}
whence $V(I) \subset V(f_2)$.
Let $\eta$ be any point of the open part of $V(f_2)$ where $a_3\ne 0$ and $f_3\ne 0$.
Then we have $\Delta = 0$ and $a_1x+a_4\ne 0$ at any point of $V(P,q)$ by Lemma \ref{lem:meaning_f_3}.
This means that $C$ at $\eta$ is singular, whence $V(f_2) \subset V(I)$.
\end{proof}

By this proposition, we can prove Main Theorem \ref{MainTheorem1}.
\begin{proof}[Proof of the former assertion in Main Theorem \ref{MainTheorem1}.]
Let \[
\pi: {\mathcal C} \to S
\]
be the family of our curves (i.e., $S$ is the space of $A,B,a_0,\ldots, a_5$).
Let $Z$ be the subscheme consisting of points $P$ on ${\mathcal C}$ such that $\pi^{-1}(\pi(P))$ is singular at $P$.
By the above proposition, $\pi(Z)$ is dense in $V(f_2)$.
Note that $Z$ is closed, since
the condition that $\pi^{-1}(\pi(P))$ is singular at $P$ for $P \in {\mathcal C}$ is described by algebraic equations.
Then since $\pi$ is projective,
$\pi(Z)$ is closed and therefore $\pi(Z) = V(f_2)$.
\end{proof}

\begin{proof}[Proof of the {latter} assertion in Main Theorem \ref{MainTheorem1}.]
We fix an elliptic curve $E=V(P)$ in the following.
We consider the natural action of the symmetric group $\mathfrak{S}_6$ of degree $6$
on $E^6$.
Put
{
    \begin{align*}
    B_0 &:= \left\{ (P_1,\ldots,P_6)\in E^6 \mid P_1 = P_2, \sum\nolimits_{i=1}^{6}{P_i} = O \right\},\\
    B &:= \left\{(P_1,\ldots,P_6)\in E^6 \mid P_i=P_j \text{ for some } i\ne j, \text{and} \sum\nolimits_{i=1}^{6}{P_i} = O\right\}.
    \end{align*}
}
It is clear that $B_0$ is geometrically irreducible and $B$ is the union of $\sigma B_0$ for $\sigma \in \mathfrak{S}_6$. The moduli of our nonsingular curves of the form $V(P,Q)$ is given by $(E^6\smallsetminus B)/\mathfrak{S}_6$ and
the loci $V(f_2)$ consisting of singular curves $V(P,Q)$ is described as $B/\mathfrak{S}_6$.
By the argument above, $B/\mathfrak{S}_6$ is an image of the geometrically irreducible variety $B_0$.
This implies that $V(f_2)$ is geometrically irreducible.
Considering these as a family where $P$ varies, clearly $V(f_2)$ is still geometrically irreducible. It remains to show that the ideal $(f_2)$ is radical. This follows from the fact that the lowest term 
of $f_2$ with respect grevlex with
\[
    A \succ B \succ a_0 \succ a_1 \succ a_2 \succ a_3 \succ a_4\succ a_5
\]
is $-6^6a_3^5a_5^7$ (we checked this by using the code in \cite{webpage}), which is not an $n$-th power of any monomial for any $n\ge 2$.
\end{proof}

\begin{Rem}
It is possible to confirm that $f_2$ is irreducible over ${\mathbb Q}$,
by using the built-in function \texttt{IsPrime} (or \texttt{Factorization}) of MAGMA \cite{MAGMA}, see \cite{webpage}.
\end{Rem}

{
\begin{Cor}
With the notion of Proposition \ref{prop:SingularLocus} (1), we have
$I = (f_2)$.
\end{Cor}
\begin{proof}
By Proposition \ref{prop:SingularLocus} (1), its proof and the latter assertion of Main Theorem \ref{MainTheorem1}, we have $(f_2) \subset I \subset \sqrt{I} = \sqrt{(f_2)} = (f_2)$.
\end{proof}
}

\section{Proof of Main Theorem \ref{MainTheorem2}}\label{sec:ProofMainTheorem2}

\label{sec:computational_results}
In this section, we give a proof of Main Theorem \ref{MainTheorem2}, a computational result for enumerating superspecial UDCEC's for small characteristic $p \leq 23$:

In Section 4.1, we give an explicit algorithm (Algorithm \ref{listDCEC}) for listing superspecial DCEC's. In Section 4.2, we explain how to enumerate superspecial UDCEC's and give the computational results (Table \ref{table:DCEC}) by executing our method. As supplements, we write down defining polynomials of superspecial UDCEC's for $p=11,13$ in Sections 4.4-4.5. Throughout this section, let $p$ be a prime integer $p \geq 5$ and $K$ the algebraic closure $\overline{{\mathbb F}_p}$ of the finite field $\mathbb{F}_{p}$.

\subsection{How to list superspecial DCEC's}\label{sebsec:Enumeration_SSP-DCEC}
In this subsection, we give an explicit algorithm (Algorithm \ref{listDCEC}) to list up superspecial DCEC's.

\begin{algorithm}
\begin{algorithmic}[1]
\renewcommand{\algorithmicrequire}{\textbf{Input:}}
\renewcommand{\algorithmicensure}{\textbf{Output:}}
\caption{Listing superspecial DCEC's}
\label{listDCEC}
\Require A prime integer $p > 3$
\Ensure The list $\mathcal{L}$ of superspecial DCEC's in characteristic $p$
\State Let $\mathcal{L} := \emptyset$
\State Let $\mathcal{E}$ be the list of the isomorphism classes of supersingular elliptic curves {in characteristic $p$}
\For{$E = V(P) \in \mathcal{E}$}
\State Let $Q'$ be the polynomial as in (\ref{PQ}) with indeterminate $a'_0,\ldots,a'_5$
\State Compute the polynomial $f_2$ given in Main Theorem \ref{MainTheorem1}
\State Compute the {Hasse-Witt} matrix $M = (c_{ij})_{i,j}$ of $C = V(P,Q')$ by using Algorithm \ref{carier-manin:algrm}
\State Let $I \subset \mathbb{F}_{p^2}[a'_0,\ldots,a'_5,s']$ be the ideal generated by
\[
    c_{ij}~(1 \leq i,j \leq 4)\quad {\rm and} \quad s'f_2 - 1,
\]
where $s'$ is a ``dummy" indeterminate
\State {Compute the extension field $L$ of $\mathbb{F}_{p^2}$ splitting all the single variable $g\in I$}
\State {Compute the set $V:=V(I)$ considered over $L$} 
\For{$(a_0,\ldots,a_5,s) \in V$}
\State Let $Q$ be a polynomial as in (\ref{PQ}) {determined by} $a_0,\ldots,a_5$
\State Add a DCEC $C = V(P,Q)$ to the list $\mathcal{L}$
\EndFor
\EndFor\\
\Return $\mathcal{L}$
\end{algorithmic}
\end{algorithm}

In the following, we give a brief description of Algorithm \ref{listDCEC}. First, we construct a set ${\mathcal E}$ of representatives of the isomorphism classes of supersingular elliptic curves such that any element of ${\mathcal E}$ is of the form $V(P)$ with {$P=y^2 z-(x^3+Ax^2z+Bz^3)$}. Next, for each $V(P)\in {\mathcal E}$, we consider the coefficients $a'_0,\ldots,a'_5$ of defining equations \eqref{PQ} of DCEC $C$ as indeterminates,
and make the simultaneous equations representing that the {Hasse-Witt} matrix of $C$ is the zero matrix and that $C$ is nonsingular. Let $I$ be the ideal of the simultaneous equations.
Finally, we solve the simultaneous equations. We note that Main Theorem \ref{MainTheorem1} ensures that the DCEC's obtained by Algorithm \ref{listDCEC} are all nonsingular.

In our implementation with computer algebra system MAGMA\,\cite{MAGMA}, we use the built-in function {\tt Variety} to solve the simultaneous equations in {Step 9.} However, the function {\tt Variety} outputs the set of rational points over the base field (a finite field), that is, we cannot get all solutions over the algebraic closure of it. For this, we list every polynomial in one variable (with respect to each indeterminate) belonging to $I$, by finding an elimination ideal, and consider the extension field $L$ splitting all the polynomials, and apply {\tt Variety} to the ideal $I$ considered over $L$. 


\begin{Rem}\label{radical-remark}
Let $I_M$ be the ideal of $K[a_0,\ldots, a_5,1/f_2]$ generated by entries of the {Hasse-Witt matrix} $M$. In the computation for Main Theorem \ref{MainTheorem2}, we verify that
the number of points on $V(I_M)$
is equal to the dimension of $K$-vector space
\[
    K[a_0,\ldots, a_5, 1/f_2]/ I_M,
\]
and therefore $I_M$
is radical for all primes $p$ with $5 \leq p \leq 23$ (see \cite{webpage} for details).
This can be also checked by using the built-in function {\tt IsRadical} of MAGMA\,\cite{MAGMA}.
{It would be natural to expect that $I_M$ is radical for all $p \geq 5$ (Conjecture \ref{radical}).}
\end{Rem}

\subsection{Whole algorithm and computational result}\label{subsec:WholeAlgorihm}
In this subsection, we explain how to enumerate superspecial DCEC's and UDCEC's, and give the computational result obtained by executing it. Here are two algorithms for enumerating superspecial DCEC's (Algorithm \ref{enumerateUDCEC}) and UDCEC's (Algorithms \ref{enumerateUDCEC} and \ref{enumerateUDCEC2}):

\begin{algorithm}
\begin{algorithmic}[1]
\renewcommand{\algorithmicrequire}{\textbf{Input:}}
\renewcommand{\algorithmicensure}{\textbf{Output:}}
\caption{Enumerating superspecial (U)DCEC's}
\label{enumerateUDCEC}
\Require A prime integer $p > 3$
\Ensure The list $\mathcal{D}$ (resp. $\mathcal{U}$) of all isomorphism classes (over $\overline{{\mathbb F}_p}$) of superspecial DCEC's (resp. UDCEC's) in characteristic $p$
\State Make the list $\mathcal{L}$ of superspecial DCEC's in characteristic $p$ using Algorithm \ref{listDCEC}
\State Let $\mathcal{D} \leftarrow \emptyset,\ \mathcal{U} \leftarrow \emptyset$
\ForAll{$C_1 \in \mathcal{L}$}
    \If {there does not exist $C_2 \in \mathcal{D}$ such that $C_1$ is isomorphic to $C_2$}
\State Add $C_1$ to the list $\mathcal{D}$
\EndIf
\EndFor
\ForAll{$C \in \mathcal{D}$}
    \If {the output of Algorithm \ref{checkUDCEC} for $C$ is True}
        \State Add $C$ to the list $\mathcal{U}$
    \EndIf
\EndFor\\
\Return $\mathcal{U}$ and $\mathcal{D}$
\end{algorithmic}
\end{algorithm}

\begin{algorithm}
\begin{algorithmic}[1]
\renewcommand{\algorithmicrequire}{\textbf{Input:}}
\renewcommand{\algorithmicensure}{\textbf{Output:}}
\caption{Enumerating superspecial UDCEC's}
\label{enumerateUDCEC2}
\Require A prime integer $p > 3$
\Ensure The list $\mathcal{U}$ of all isomorphism classes (over $\overline{{\mathbb F}_p}$) of superspecial UDCEC's in characteristic $p$
\State Make the list $\mathcal{L}$ of superspecial DCEC's in characteristic $p$ using Algorithm \ref{listDCEC}
\State Let $\mathcal{L}' \leftarrow \emptyset,\ \mathcal{U} \leftarrow \emptyset$
\ForAll{$C \in \mathcal{D}$}
    \If {the output of Algorithm \ref{checkUDCEC} for $C$ is True}
        \State Add $C$ to the list $\mathcal{L}'$
    \EndIf
\EndFor
\ForAll{$C_1 \in \mathcal{L}'$}
    \If {there does not exist $C_2 \in \mathcal{U}$ such that $C_1$ is isomorphic to $C_2$}
\State Add $C_1$ to the list $\mathcal{U}$
\EndIf
\EndFor\\
\Return $\mathcal{U}$
\end{algorithmic}
\end{algorithm}

The main difference between these two algorithms is the order of the isomorphism tests and the step determining whether a DCEC is a UDCEC.

\begin{Rem}
In Algorithm \ref{enumerateUDCEC2}, since all the curves in $\mathcal{L}'$ are UDCEC's, we can use Algorithm \ref{isomtest} to test whether $C_1 \cong C_2$ for all $C_1,C_2 \in \mathcal{L}'$. On the other hand, in Algorithm \ref{enumerateUDCEC}, we cannot use Algorithm \ref{isomtest}. Instead, we use the built-in function {\tt IsIsomorphic} in MAGMA or \cite[Sections 4.1, 4.2 and 6.1]{KH16}.
\end{Rem}

{We chose to implement} Algorithm \ref{enumerateUDCEC} in this paper. Using Algorithm \ref{enumerateUDCEC}, it is possible to compute not only the number of isomorphism classes of superspecial UDCEC's, but also the {number} of superspecial DCEC's. Of these, the number of isomorphism classes of superspecial Howe curves for $p < 20000$ is known by Kudo-Harashita-Howe \cite{KHH}. {Using Algorithm \ref{enumerateUDCEC} lefts us verify that our results are consistent.}

\begin{table}[htbp]
  \begin{tabular}{cccccccc}
   \toprule             

   char.  & s.sing $j$-inv. & \multicolumn{4}{c}{the number of isom.\ classes}  \\
   \cmidrule(lr){3-6}   
      & & all & UDCEC  & Howe & others  \\
   \midrule             
    $p=5$ & 0          & 1 & 0 & 1 & 0  \\
    $p=7$ & 1728       & 0 & 0 & 0 & 0  \\
    $p=11$ & 0,1728    & 5 & 1 & 4 & 0  \\
    $p=13$ & 5         & 5 & 2 & 3 & 0  \\
    $p=17$ & 0,8       & 27 & 13 & 10 & 4   \\
    $p=19$ & 7,1728    & 20 & 16 & 4 & 0   \\
    $p=23$ & 0,19,1728 & 95 & 53 & 33 & 9  \\   
   \bottomrule          
  \end{tabular}
   \caption{The number of isomorphism classes of all superspecial DCEC's in characteristic $p \leq 23$
   }
 \label{table:DCEC}
\end{table}

We implemented and executed Algorithm \ref{enumerateUDCEC} {using the} Computational Algebra System MAGMA \cite{MAGMA}, with the version V2.26-12. Our codes and log files are available at the website \cite{webpage}. Table \ref{table:DCEC} {summarizes} our computational results for characteristic $p \leq 23$. Our implementations over MAGMA were conducted on a PC equipped with an Ubuntu 22.04.3 LTS OS, utilizing a 4.00GHz $\times 8$ quad-core CPU (Intel(R) Xeon(R) W-2125) and 128GB of RAM. 
On the computation, the least time costs for $p=5$ of all characteristics above, and the time is 9.97 seconds. The most time costs for $p=23$, and the time is 108823 seconds, approximately 30 hours and 13 minutes.

\begin{Rem}
In Table \ref{table:DCEC} the above, we classify all DCEC's $C = V(P,Q)$ into the following three types:
\begin{itemize}
\item If $C$ has a unique elliptic quotient of degree 2, then $C$ is a UDCEC.
\item If $C$ has another elliptic quotient of degree 2 and the associated two involutions in ${\rm Aut}(C)$ commute (in other words,
the automorphism group of $C$ contains $V_4$), then $C$ is a Howe curve, see Proposition \ref{HoweV4}.
\item Otherwise $C$ is classified as “others".
In case, there are at least two involutions in ${\rm Aut}(C)$
{that do not} commute. Note that the two involutions generate a dihedral group of order $\ge 6$.
\end{itemize}
Note that the number of isomorphism classes of superspecial Howe curves for characteristic $p < 20000$ was already given by \cite{KHH}, and Table \ref{table:DCEC} is consistent with their result.
\end{Rem}

From the above discussions, the proof of Main Theorem \ref{MainTheorem2} is done. In the rest of this section, we describe these results more precisely, including explicit defining equations of them.

\subsection{Remarks for $p \leq 7$}
In characteristics $p \leq 7$, our computation shows that there is no superspecial UDCEC, which is not a new result. In the following, let us explain this briefly:
\begin{itemize}
\item Ekedahl \cite[Theorem 1.1]{Ekedahl} showed that if there exists a superspecial curve of genus $g$ in characteristic $p$, then $2g \leq p^2 - p$ holds. This implies that there does not a superspecial curve of genus $4$ in characteristic $p=3$.
\item Fuhrmann-Garcia-Torres \cite[Theorem 3.1]{FGT} showed that there exists a unique superspecial non-hyperelliptic genus-4 curve in characteristic $p=5$, which turned out to be a Howe curve.
\item Kudo-Harashita \cite{KH16} showed that there does not exist a superspecial curve of genus $4$ in characteristic $p=7$.
\end{itemize}
Therefore, Table \ref{table:DCEC} for characteristic $p \leq 7$ is consistent with their results.

\subsection{On results for $p = 11$}
In characteristic $p = 11$, there exist two supersingular $j$-invariants: $j=0$ or $1728$. For $j \in \{0,1728\}$, we let $E$ be a supersingular elliptic curve with $j(E) = j$.
\begin{enumerate}
\item The case of $j = 0$.\\
In case, $E$ is given by 
$P := y^2z - (x^3 + xz^2) = 0$. 
Moreover, there does not exist any superspecial UDCEC $C = V(P,Q)$.
\item The case of $j = 1728$.\\
In case, $E$ is given by 
$P := y^2z - (x^3 + z^3)=0$. Moreover, there exists a unique superspecial UDCEC $C = V(P,Q)$ whose explicit equation is given as
\[
    C: \ \left\{
    \begin{array}{l}
        y^2z = x^3 + z^3,\\
        w^2 = xy.
    \end{array}
    \right.
\]
This curve $C$ is isomorphic to $V(Q^{({\rm Dege})},{P_6}^{({\rm alc)}})$ in \cite[Proposition 5.2.5]{KH17} with
\[
    \begin{array}{l}
        Q^{({\rm Dege})}:=2YW + Z^2 = 0,\\
        {P_6}^{({\rm alc)}}:=X^3 + XW^2 + Y^3 = 0,
    \end{array}
\]
by setting $X = z,\ Y = x,\ Z = \sqrt{-2i}w$ and $W = iy$ with $i^2 = -1$.
\end{enumerate}

\subsection{On results for $p = 13$}
In characteristic $p = 13$, there exists only one supersingular $j$-invariant: $j = 5$. In case, $E$ is given by the equation
\[
    E: y^2z = x^3 + 5xz^2 + {\sqrt{-2}}z^3 
\]
where $\hspace{-0.3mm}{\sqrt{-2}} \in\mathbb{F}_{13^2}$ is one of roots of the quadratic equation $s^2 = {-2}$. Set $P := y^2z - (x^3 + 5xz^2 + {\sqrt{-2}}z^3)$. Then, there exist {two isomorphism classes of} superspecial UDCEC's $C = V(P,Q)$. One of them is given as
\[
    C_1: \left\{
    \begin{array}{l}
        y^2z = x^3 + 5xz^2 + {\sqrt{-2}}z^3,\\
        w^2 = x^2 + \alpha xy + 6{\sqrt{-2}}xz + \bigl(8{\sqrt{-2}}+3\bigr)y^2 + 3\alpha{\sqrt{-2}} y z + \bigl(8{\sqrt{-2}}- 6\bigr)z^2,
    \end{array}
    \right.
\]
where $\alpha$ denotes a root $\in\mathbb{F}_{13^4}$ of the quadratic equation $t^2 = 6{\sqrt{-2}}-1$. Another curve is given as
\[
    C_2: \left\{
    \begin{array}{l}
        y^2z = x^3 + 5xz^2 + {\sqrt{-2}}z^3,\\
        w^2 = x^2 + \beta xy + 6{\sqrt{-2}}xz + \bigl(8{\sqrt{-2}}-3\bigr)y^2 + 3\beta{\sqrt{-2}} y z - \bigl(8{\sqrt{-2}}+6\bigr)z^2,
    \end{array}
    \right.
\]
where $\beta$ denotes a root $\in\mathbb{F}_{13^4}$ of the quadratic equation $t^2 = 6{\sqrt{-2}}+1$.
{\begin{Rem}
There are two choices for $\alpha$ (resp.~$\beta$) as the root of a quadratic equation, but the two values of $\alpha$ (resp.~$\beta$) give isomorphic curves over $\mathbb{F}_{13^4}$.
\end{Rem}}

\section{Conclusions and future work}\label{sec:conclusions}
We studied the defining equation of a DCEC (genus-4 double cover of an elliptic curve) $C$, and we gave a necessary and sufficient condition for $C$ to be nonsingular,
which is described as $f_1\ne 0$ and $f_2\ne 0$
with discriminant $f_1$ of the elliptic curve.
Note that $f_2$ is a very long polynomial (with 3381 terms in characteristic $0$), but
we were able to show that $f_2$ is an irreducible polynomial. These results are stated in {Main Theorem 1, which is a theoretical result} working over arbitrary field of characteristic $\ne 2,3$.

As an application, we had a computational result on
the enumeration of superspecial DCEC's in small characteristic $p \leq 23$. 
In Main Theorem \ref{MainTheorem2}, we focused on the case that there is a unique morphism $C$ to an elliptic curve of degree two.
Such a curve is called UDCEC.
Other cases (e.g. the case of Howe curves) can be discussed by using other defining equations, since they are special cases.
We proposed an algorithm to list superspecial DCEC's, an algorithm to select UDCEC's among them,
and an isomophism-test algorthm for (U)DCEC's.
We implemented these algorithms {in the} computer algebra system MAGMA, and succeeded in enumerating isomorphism classes of superspecial UDCEC's in characteristic $p\le 23$.

{In the computarion for Main Theorem \ref{MainTheorem2}, we have found the following interesting conjecture, as stated in Remark \ref{radical-remark}:
\begin{Conj}\label{radical}
With the same notations in Section \ref{subsec:Preliminaries-DCEC}, let $I_M$ be the ideal of $K[a_0,\ldots,a_5,1/f_2]$ generated by entries of the Hasse-Witt matrix of $C$. Then, the ideal $I_M$ is radical for all $p \geq 5$.
\end{Conj}
\noindent This conjecture can be regarded as an analogue of Igusa's result \cite{Igusa}: the entry $H_p(t)$ of a {Hasse-Witt} matrix ($1 \times 1$ matrix) of the elliptic curve $y^2=x(x-1)(x-t)$ is a separable polynomial.
For an analogous result in the genus-$2$ case, see \cite[Prop.\ 1.14]{IKO} 
in some special cases and
\cite{HY} for genus-$2$ curves in Rosenhain form.}

{The following would be in the interesting future directions}:
\begin{enumerate}
    \item {An analogous result for Main Theorems \ref{MainTheorem1}} for characteristic $p=2,3$. 
    \item To find invariants of DCEC's (or UDCEC's) would be important to give more effective isomorphism-test algorithm.
    \item Is there theoretical approach toward the enumeration of superspecial DCEC's (or UDCEC's)? 
  \end{enumerate}
  
    A possible approach for the third item would be to find {a} relation to genus-3 curves. Indeed, the moduli space of DCEC's and that of genus-3 curves {are both} dimension $6$. For a DCEC $C$, let $C \to E$ be a double cover with elliptic curve $E$. The Prym variety $A:={\rm Ker}({\rm Jac}(C) \to E)$ is a polarized abelian variety of dimension $3$.
    Ikeda \cite{Ikeda} proved
    that this correspondence $C \mapsto A$ is injective. Of course, if $C$ is superspecial, then $A$ is superspecial. (Unfortunately the converse seems to be false in general.) 
    We also pay attention to another correspondence by Del Centina and Recillas \cite{CR}
    (also see Dolgachev \cite[Section 3]{Dolgachev}), which says that the birational isomorphism from 
    the moduli space ${\mathcal R}_3$ of genus-$3$ curves with a non-trivial $2$-tortion divisor class
    to the moduli space ${\mathcal R}_4^{\text{be}}$ of genus-$4$ curves with a bi-elliptic involution $\sigma$
    and a $\sigma$-invariant $2$-torsion divisor class.
  %
  With these facts, it is natural to expect
  that the number of isomorphism classes of superspecial Prym varieties $A$ as above
  is proportional to the number of isomorphism classes of superspecial genus-$3$ curves as $p\to\infty$.
    If we assume
    that there is no relation between the superspeciality of $A$ and that of $E$, then
    the facts that the number of isomorphism classes of superspecial genus-$3$ curves is approximately $p^6/1451520$ (\cite{Hashimoto} and \cite{Brock}) and the {probability} for elliptic curves $E$ over ${\mathbb F}_{p^2}$ to be supersingular is approximately $(12p)^{-1}$ as $p\to\infty$ lead us to propose the second expectation of the following:
   \begin{Expect}\label{expect:AsymptoticBehavior}
   \begin{enumerate}
   \item[(1)] There exists a superspecial UDCEC in every characteristic $p \ge 11$.
    \item[(2)]  The number of isomorphism classes of superspecial (U)DCEC's would be $O(p^5)$ 
    as $p\to\infty$.
   \end{enumerate}
   \end{Expect}
Here is a table for $p\le 23$.
We exclude the contribution from Howe curves and others,
since it will likely diminish, compared with that of UDCEC's, as $p$ increases 
   (approximately $p^3/1152$ for the Howe curves as expected in \cite{KHH}).
As $p$ is small, it seems insufficient for estimating the proportional constant.
Here, we tentatively set it to $1/120000$.
      \begin{table}[htbp]
  \begin{tabular}{|c|c|c|c|c|c|c|c|}
   \hline
   $p$  & $5$ & $7$ & $11$ & $13$ & $17$ & $19$ & $23$  \\
   \hline
  \# of ssp. UDCEC's & 0 & 0 & 1 & 2 & 13 & 16 & 53\\
  \hline
  $p^5/120000$ & 0.0 &  0.1 & 1.3 & 3.1 & 11.8 & 20.6 & 53.6\\
  \hline
  \end{tabular}
   \caption{The number of isomorphism classes of all superspecial UDCEC's vs. $p^5/120000$
   }
\end{table}


\begin{thebibliography}{99}
%
%

\bibitem{MAGMA}
Bosma, W., Cannon, J.\ and Playoust, C.: \textit{The MAGMA algebra system. I. The user language}, J. Symb. Comput.\,{\bf 24}, No.\,3-4, 235--265 (1997).


\bibitem{Brock} Brock, B. W.:
\textit{Superspecial curves of genera two and three},
Thesis (Ph.D.)-Princeton University (1993).


\bibitem{CR}
Del Centina, A.\ and Recillas, S.:
\textit{On a property of the Kummer variety and a relation between two moduli spaces of curves},
Algebraic geometry and complex analysis, Proc. Workshop, P\'atzcuaro/Mexico 1987, Lect. Notes Math. 1414, 28-50 (1989).







\bibitem{Deuring}
Deuring, M.:
\textit{Die Typen der Multiplikatorenringe elliptischer Funktionenk\"orper},
Abh. Math. Sem. Univ. Hamburg {\bf 14}, No.\,1, 197--272 (1941). 

\bibitem{Dolgachev}
Dolgachev, I.\ V.:
\textit{Rationality of ${\mathcal R}_2$ and ${\mathcal R}_3$},
 Pure Appl.\ Math.\ Q.~4, No.~2, 501--508 (2008).
 



\bibitem{Ekedahl}
Ekedahl, T.: \textit{On supersingular curves and abelian varieties},
Math. Scand.\,{\bf 60}, 151--178 (1987).



\bibitem{FGT} Fuhrmann, R., Garcia, A., Torres, F.:
{\it On maximal curves},
J. Number Theory {\bf 67}, No. 1, 29--51 (1997).







\bibitem{HY}
Harashita, S.\ and Yamamoto, Y.:
{\it The multiplicity-one theorem for the superspeciality of curves of genus two}, arXiv:2409.13212 [math.AG]


\bibitem{Har}
Hartshorne, R.: {\it Algebraic Geometry}, GTM {\bf 52}, Springer-Verlag (1977).




\bibitem{Hashimoto}
Hashimoto, K.:
\textit{Class numbers of positive definite ternary quaternion hermitian forms}, Proceed. Japan Acad. {\bf 59} Ser. A, 490--493 (1983).



\bibitem{Howe}
Howe, E. W.:
{\it Quickly constructing curves of genus 4 with many points},
Contemp. Math.{\bf 663}, 149--173 (2016).



\bibitem{Ibukiyama}
Ibukiyama, T.:
{\it On rational points of curves of genus $3$ over finite fields},
Tohoku Math. J. {\bf 45}, 311-329 (1993),.


\bibitem{IKO}
Ibukiyama, T., Katsura, T. and Oort, F.:
{\it Supersingular curves of genus two and class numbers}, 
Compositio Math. {\bf 57}, No. 2, 127--152 (1986). 

\bibitem{Igusa}
Igusa, J.:
\textit{Class number of a definite quaternion with prime discriminant},
Proc. Nat. Acad. Sci. U.S.A. {\bf 44}, 312--314 (1958).

\bibitem{Ikeda}
Ikeda, A.:
\textit{Global Prym–Torelli theorem for double coverings of elliptic curves},
Algebraic Geometry {\bf 7} (5) (2020) 544–-560.

\bibitem{KR}
Kani, E. and Rosen M.:
\textit{Idempotent relations and factors of Jacobians},
Math. Ann. {\bf 284} (1989), 307--327. 









\bibitem{KH16}
Kudo, M.\ and Harashita, S.: \textit{Superspecial curves of genus $4$ in small characteristic}, Finite Fields and Their Applications {\bf 45}, 131--169 (2017).

\bibitem{KH17}
Kudo, M.\ and Harashita, S.: \textit{Computational approach to enumerate non-hyperelliptic superspecial curves of genus 4}, 
Tokyo Journal of Mathematics Volume {\bf 43}, Number 1, 259--278 (2020).



\bibitem{KH18}
Kudo, M.\ and Harashita, S.:
{\it Superspecial Hyperelliptic Curves of Genus $4$ over Small Finite Fields},
In: L.\ Budaghyan, F.\ Rodriguez-Henriquez (eds), Arithmetic of Finite Fields, WAIFI 2018, Lecture Notes in Computer Science, Vol.\ {\bf 11321}, pp.\ 58--73, Springer, Cham, 2018, (DOI) 10.1007/978-3-030-05153-2\_3.





\bibitem{KHH}
Kudo, M., Harashita, S. and Howe, E.\ W.:
{\it Algorithms to enumerate superspecial Howe curves of genus $4$},
Proceedings of Fourteenth Algorithmic Number Theory Symposium (ANTS-XIV), MSP, 2020. Open Book Series Vol. {\bf 4} (2020), No. 1, 301-316 

\bibitem{KHS}
Kudo, M., Harashita and S., Senda, H.:
{\it The existence of supersingular curves of genus 4 in arbitrary characteristic}, Res.\ in Number Theory {\bf 6}, Article num.: 44 (2020).










\bibitem{Nygaard}
Nygaard, N.\ O.:
{\it Slopes of powers of Frobenius on crystalline cohomology}.
Ann.\ Sci.\ \'{E}c.\ Norm.\ Sup\'{e}r.\ (4) 14, 369--401 (1981).


\bibitem{KNT}
Kudo, M., Nakagawa T.\ and Tsuyoshi T.:
{\it Efficient search for superspecial hyperelliptic curves of genus four with automorphism group containing $C_6$.}
Math. Comput. Sci.\,{\bf 17}, No. 3--4, Paper No. 21 (2023).

\bibitem{Lachaud}
Lachaud, G.:
{\it Sommes d’Eisenstein et nombre de points de certaines courbes alg\'ebriques sur les corps finis.}
C. R. Acad. Sci., Paris, S\'er. I {\bf 305}, 729--732 (1987).

\bibitem{OKH23}
Ohashi, R., Kudo, M. and Harashita, S.:
{\it Fast enumeration of superspecial hyperelliptic curves of genus $4$ with automorphism group $V_4$,}
In: Mesnager, S., Zhou, Z. (eds) Arithmetic of Finite Fields. WAIFI 2022. Lecture Notes in Computer Science, vol. {\bf 13638}. pp.107-124, Springer, Cham.

\bibitem{OK23}
Ohashi, R. and Kudo, M.: {\it Computing superspecial hyperelliptic curves of genus 4 with automorphism group properly
containing the Klein 4-group}, Journal of Computational Algebra {\bf 11} (2024), Article ID: 100020.

\bibitem{Oort-Hyp}
Oort, F.:
\textit{Hyperelliptic supersingular curves},
Progress in Mathematics, {\bf 89}, pp.~247--284, Birkh\"auser, Boston (1991).

\bibitem{OU}
Oort, F. and Ueno, K.: Principally polarized abelian varieties of dimension two or three are Jacobian varieties. Journ. Fac. Sc. Univ. Tokyo, {\bf 20}, 377-–381 (1973).






\bibitem{Serre1983}
Serre, J.-P.:
{\it Nombre des points des courbes algebrique sur ${\mathbb F}_q$},
S\'em. Th\'eor. Nombres Bordeaux (2)
1982/83, 22 (1983).

%

\bibitem{Silverman}
Silverman J. H.: \textit{The Arithmetic of Elliptic Curves}, GTM\,{\bf 106}, Springer--Verlag (1986).






%
%

\bibitem{webpage}
Codes and log files for the paper ``Superspecial genus-$4$ double covers of elliptic curves", \url{http://www.h-lab.ynu.ac.jp/SSP-genus4-DCEC.html}.

\end{thebibliography}
\end{document}